\documentclass[s5,11pt]{article}
\usepackage{amsmath}
\newtheorem{theorem}{Theorem}

\usepackage{latexsym}

\newtheorem{lemma}{Lemma}

\newtheorem{proposition}{Proposition} 
  
\author{ Tord Sj\"odin\footnote{Department of mathematics and mathematical statistics, Ume\aa\ university, S-90187 Ume\aa, Sweden; e-mail: tartansson48@gmail.com;
phone: +46-90-786 5000; Fax: -}}
\title{ The Power Problem for Generalized Gamma Convolutions (GGC) and Related Questions \footnote{Dedicated to the memory of  my colleague and friend professor emeritus Lennart Bondesson, Ume\aa\ university, Sweden}}
 \begin{document}  
 \maketitle
\begin{abstract} The class of generalized gamma convolutions (GGC) is closed with respect to change of scale, weak limits and addition and multiplication of independent random variables. Our main result confirms an old conjecture that GGC is also closed wrt $q-$th powers, $q>1$. The proof uses explicit iterative formulas for the densities of finite sums of independent gamma variables,   hyperbolically completely monotone functions (HCM) and the Laplace transform.\\[0.5em] We apply the result to sums and products of $q-th$ powers of independent GGCs, $q\geq 1$, symmetric extended GGC (symEGGC) and a new proof that $X\sim GGC$ implies $Exp(X)\sim GGC$.\end{abstract}
\paragraph{ \it AMS 2010 Subject Classsification:} Primary 60E10
Secondary 62E15
\section{Introduction}The generalized gamma convolutions (GGC) were introduced by O. Thorin \cite{T1},\cite{T2} in his study of infinite divisibility of the lognormal distribution, see also \cite{T3}. The class GGC consists of limit distributions of finite sums of independent gamma random variables (rvs) and is closed with respect to (wrt) change of scale, weak limits and sums and products of independent rvs. A comprehensive study of GGC and its relation to hyperbolically completely monotone functions (HCM) is found in Bondesson \cite{B1}, see also Steutel, van Hahn \cite{SH}, Ch. VI, \S 5 and Bondesson \cite{B2}. We use Feller \cite{F} as a general reference on probability theory. For more on the background in infinite divisibility, GGC and the pioneering work of O. Thorin, see the nice biography by Bondesson, Grandell, Peetre \cite{BGP}. A problem on a class of mixtures of gamma distributions in the same field was studied in Behme, Bondesson \cite{B3} and by the author in $\cite{S}$.\\[1em] 
Our main result (Theorem 1) confirms an old conjecture for $GGC$ going back at least to the late 1980's, that $X\sim GGC$ and $q>1$ implies that $X^q\sim GGC$ (here called the Power Problem) mentioned in Bondesson \cite{B1}, p. 97. It is known to hold in several special cases, see Bondesson \cite{B1}, Ch.6 and \cite{B2}, Sec.7. If $PF_\infty$ denotes the class of limit distributions of finite sums of exponential rvs, then $X\sim PF_\infty$ implies that $X^q\sim GGC$  for $q\geq 1$, Bondesson \cite{B1}, Theorem 6.2.7. The conjecture is then also true for sums of independent gamma rvs whose shape parameters are positive integers. A positive answer to the Power problem was conjectured in Bondesson \cite{B2}, Conjecture 1. Our main result (Theorem 1) confirms the conjecture. The proof is based and Bondessons characterization of GGC in {\cite{B1}, Theorem 5.4.1, explicit iterative formulas for the densities of finite sums of independent gamma rvs and successive substitutions.  The result is applied to a new proof that $X\sim GGC$ implies that $e^X-1\sim GGC$, Bondesson. \cite{B2}, Theorem 4, (Theorem 2), to sums  and products of powers of independent GGCs (Theorem 3) and inclusion theorems for symmetric extended GGCs, $sym EGGC$ (Theorem 4).
\\[1em]   Section 2 begins with the standard notation used in this field, a review of our set up and   three lemmas, where Lemma 2 is used in the induction step of the proof. Our main result (Theorem 1) is stated and proved in Section 3 and the applications are given in Section 4.  
 \section{Background}  This section gives the necessary background and defines the concepts needed to state and prove our theorems, c.f. Bondesson \cite{B1}. A function $f:(0,\infty )^n\rightarrow [0,\infty )$ is completely monotone (CM) if $(-1)^m D^m f \geq 0$, for  all positive integers $m$,   and a function $f:(0,\infty )\rightarrow [0,\infty )$ is  hyperbolically completely monotone (HCM) if, for every fixed $u>0$, $H(w)=f(uv)\cdot f(u/v)$ is CM wrt $w=v+v^{-1}$, see Bondesson \cite{B1}, Ch. 5. We let $Gamma(\beta ,b)$ denote a standard gamma distribution with density $f(x)= b^\beta  \cdot \Gamma (\beta )^{-1}\cdot x^{\beta -1}\cdot e^{-bx}$, $x>0$, and write $Gamma(1 ,b)=Exp(b)$ for the exponential distribution. A generalized gamma convolution (GGC) is defined as a limit distribution of finite sums $X_1+X_2+\cdots +X_n$ of independent gamma rvs $X_i\sim Gamma (\beta_i,b_i)$, $1\leq i\leq n$. Then $X\sim GGC$ if and only if the Laplace transform $\phi$ of the distribution of $X$ can be represented as
$$\phi (s)=E[e^{-sX}] =\exp \big( -as+\int \log (\frac {t}{t+s})\, U(dt)\big),\, s\geq 0,$$
where $a\geq 0$ (called the left extremity) and $U(dt)$ is a nonnegative measure on $(0,\infty )$, with finite mass on compact subsets of $(0,\infty )$, such that $\int _0^1|\log t| \, U(dt)<\infty $ and $\int _1^\infty t^{-1}\, U(dt)<\infty$, Bondesson \cite{B1} Ch. 3.
  We use the following well-known characterization of GGC, see Bondesson \cite{B1}, Theorem 5.4.1.
 \begin{proposition} A function $\phi (s)$ defined on $(0,\infty )$ such that $\phi (0+)=1$ is the Laplace transform of a GGC if and only if $\phi$ is HCM.\end{proposition}
 We begin our analysis of $q-$th powers of $GGC$ and the proof that $X\sim GGC$ implies $X^q\sim GGC$ by considering finite sums $S_n=X_1+X_2+\cdots +X_n$ of independent gamma rvs, $X_i\sim Gamma(\beta_i,b_i)$, $1\leq i\leq n$. We recall that if all $\beta_i=1$, $1\leq i\leq n$, then $S_n$ is a sum of independent exponentially distributed rvs and $X^q\sim GGC$ by Bondesson \cite{B1}, p. 96. The same conclusion follows if all $\beta_i$ are positive integers, since then each $X_i$ is a sum of $\beta_i$ independent $Exp(b_i)$ rvs and we are back in the first case. Since $GGC$ is closed wrt weak limits it is no loss of generality to assume that each $\beta_i$ is a rational number $\beta_i=p_i/N$, for some positive integer $p_i$ and a common denominator $N\geq 2$, $1\leq i\leq n.$ Then each $X_i$ is the sum of $p_i$ independent rvs $X_{i,j}\sim Gamma(1/N,b_i)$, $1\leq i\leq n$, $1\leq j\leq p_i$. The Laplace transform $\phi_{S_n}$ of $S_n$ becomes
 $$\phi_{S_n}(s)=\prod\limits_{i=1}^n \big(\frac {b_i}{s+b_i}\big)^{p_i/N}=\prod\limits_{i=1}^n\bigg( \prod\limits_{j=1}^{p_i}\big(\frac {b_i}{s+b_i}\big)^{1/N}\bigg)$$
   and $S_n$ is a finite sum of independent gamma distributed rvs with form parameter {1/N}, by the uniqueness of the Laplace transform, Lemma 3. It is thus no loss of generality to assume that $\beta_i=\beta$, $1\leq i\leq n$, for some $ \beta >0$.   \\[0,5em] We start with the formulas for the density of the sums $S_n$ above in their most general form and specialize to the case $\beta_i=\beta$, $1\leq i\leq n$, later. When $n=2$, a  direct calculation gives
 \begin{equation}f_{S_2}(x)=\frac{b_1^{\beta_1}\cdot b_2^{\beta_2}}{\Gamma (\beta_1)\cdot \Gamma (\beta_ 2)  }\cdot x^{ \beta_1+\beta_2 -1}\cdot   \int\limits _0^1 e^{-x\cdot (b_1(1-u)+b_2u)}\cdot (1-u)^{\beta_1-1}\cdot u^{\beta_2 -1}\, du\end{equation}
 and for a general $n$ we use Akkouchi \cite{A}, Theorem 1 to get
   \begin{equation}f_{S_n}(x)  =D_n\cdot x^{\beta_1+\beta_2+\cdots +\beta_n -1}\cdot   \int\limits _0^1\cdots \int\limits _0^1e^{-x\cdot C_n({\bf u})}\cdot B_n({\bf u})\, du_1du_2\cdots  du_{n-1},\end{equation}
   where ${\bf u} =(u_1,u_2,\dots ,u_{n-1})$ and
   $$C_n({\bf u})=      b_1(1-u_1) +b_2 u_1(1- u_2)+\cdots +b_{n-2}\cdot u_1 u_2\cdots u_{n-3}\cdot (1-u_{n-2})+  $$ $$ +b_{n-1}\cdot u_1u_2\cdots u_{n-2}\cdot (1-u_{n-1}) +b_n \cdot u_1u_2\cdots u_{n-1} ,$$
   We note that the sum of last two terms in $C_n({\bf u})$ simplifies to
   \begin{equation}u_1u_2\cdots u_{n-2}\cdot \big( b_{n-1}\cdot (1-u_{n-1})+b_n\cdot u_{n-1}\big), \end{equation}
   which is used in the proof of Theorem 1 for $n\geq 3$. Further
   $$ B_n({\bf u})=\frac{\Gamma (\beta_1+\beta_2+\cdots +\beta_n)}{\Gamma (\beta_1) \Gamma (\beta_2)\cdots \Gamma (\beta_n)}\cdot \prod\limits _{j=1}^{n-1}u_j ^ {\beta_1+\beta_2+\cdots + \beta_j-1}\cdot (1-u_j)^{\beta_{j +1}-1}$$
   \textrm{and}
   $$ D_n=\frac{b_1^{\beta_1} \cdot b_2^{\beta_2} \cdots b_n^{\beta_n}}{\Gamma (\beta_1+\beta_2+\cdots +\beta_n )},$$
     for all ${\bf u}$. 
   \begin{lemma} [Feller \cite {F}, Criterium 2, p. 441] Let $f:(0,\infty )\rightarrow (0,\infty )$ be CM and assume that $g:(0,\infty )\rightarrow (0,\infty )$ has a CM derivative. Then $f\circ g$ is CM.\end{lemma}
   The next lemma is used in the induction step of the proof of Theorem 1. 

\begin{lemma} Let $b_1,b_2, \beta_1, \beta_2, A$ and $B$ be positive numbers and $0<\alpha <1$. Then
   \begin{equation}L=\int\limits_0^1\int\limits_0^1 e^{-  E}\cdot  g(u,v)   \, du\, dv,\end{equation}
   where $g=g(u,v)=\big( (1-u)(1-v)\big)^{\beta_1-1}\cdot (uv)^{\beta_2-1}$ and
   $$E= b_1\big((1-u)Ay^\alpha +(1-v)By^{-\alpha}\big)+b_2(uAy^\alpha +By^{-\alpha})$$
   is CM wrt $Ay^\alpha+By^{-\alpha}$, for $y>0$.
   \end{lemma}
   The proof of Lemma 2 for $n=2$ is contained in the proof of Theorem 1. The general case is proved at the end of the next section.    \begin{lemma}[Feller \cite{F}, Chap. XIII.I,Theorem 1, p. 408] (Uniqueness.)  Distinct probability distributions have distinct Laplace transforms.\end{lemma}
  \section{Powers  }
In this section we state and prove our main result that GGC is closed wrt taking $q-$th powers, $q>1$.\begin{theorem} Let $q>1$ and assume that $X\sim GGC$. Then $X^q\sim GGC.$
\end{theorem}
  {\it Proof.} As noted above, it is enough to prove the theorem for finite sums $S_n=X_1+X_2+\cdots +X_n $ of independent gamma variables, where $X_i\sim Gamma(\beta ,b_i)$, $1\leq i\leq n$, for all $\beta >0$. The proof is by induction over $n$ and uses explicit formulas for the density $f_{S_n}$ of $S_n$ and the Laplace transform $\phi_{S_n^q}$ of $S_n^q$.   
 We start with the case $n=2$. Then the density of $S_2$ is given by (1) and
  $S_2^q$ has Laplace transform 
 $$\phi_{S_2^q }(s)\sim \int\limits_0^\infty  e^{-sx^q}\cdot x^{2\beta -1}\cdot \bigg(   \int\limits _0^1 e^{- x\cdot (b_1\cdot (1-u)+b_2\cdot u)}\cdot  (u\cdot (1-u))^{\beta -1}\, du\bigg)\, dx.$$
 We will use Proposition 1 and recall the definition of the class HCM in Section 1. We compute $H_2=\phi_{S_2^q}(st)\cdot \phi_{S_2^q}(\frac{s}{t})$ as a product of two such integrals and get 
 \begin{equation}H_2\sim \int\limits_0^\infty\int\limits_0^\infty e^{- sx^q\cdot (t\cdot y+t^{-1}\cdot y^{-1} ) } \cdot  x^{2\beta -1}\cdot \end{equation}
 $$\int\limits_0^1\int\limits_0^1 e^{-x\cdot\big(   b_1\cdot ((1-u)\cdot y^\alpha +(1-v)\cdot y^{-\alpha})+b_2\cdot (u\cdot y^\alpha +v\cdot y^{-\alpha})\big) }\big( u(1-u)\cdot v(1-v)\big)^{\beta -1}\, dudv\, dx\frac{dy}{y},$$
 after a hyperbolic change of variables $x\rightarrow x\cdot y$, $y\rightarrow x/y$ and a substitution $y\rightarrow y^\alpha$, $\alpha =1/q$,
 and set out to prove that $H_2$ is CM wrt $t+t^{-1}$. We fix $s$ and $x$ and denote the inner integral by $I_2$. Then, after substitutions $u\rightarrow  \frac{1}{1+u}$ and $v\rightarrow  \frac{1}{1+v}$, $I_2=\int\int e^{-E_1}\cdot g(u,v)du dv$,  where
 $$E_1=x\cdot  \frac{b_1( u\cdot y^\alpha +v\cdot y^{-\alpha})+b_2( v\cdot y^\alpha +u\cdot y^{-\alpha})+(b_1uv+b_2)\cdot (y^\alpha +y^{-\alpha}) }{(1+u)(1+v)}$$
 and $g(u,v)=\frac{(uv)^{\beta -1}}{((1+u)(1+v))^{2\beta}}$ and the integration is over $(0,\infty )\times (0,\infty )$.
We set out to prove that $I_2$ is CM wrt $y^\alpha +y^{-\alpha}$. Without loss of generality, we assume that $y>1$ and put $y^\alpha +y^{-\alpha}=2s$. Then we get 
$$y^\alpha =s+\sqrt{s^2-1}\textrm{ and } y^{-\alpha}=s-\sqrt{s^2-1}$$ and note that also $y^{\alpha}$ is a Bernstein function wrt $y^\alpha +y^{-\alpha}$.\\[0.5em]
Next we define  
 $$E_2=  b_1( u\cdot y^\alpha +v\cdot y^{-\alpha})+b_2( v\cdot y^\alpha +u\cdot y^{-\alpha}) = (b_1u+b_2v)\cdot y^\alpha +(b_1v+b_2u)\cdot y^{-\alpha }$$ 
 and denote
 $$\Delta =(b_1u+b_2v)- (b_1v+b_2u)=(b_1-b_2)(u-v).$$
 We observe that the integral $I_2$ is unchanged if $b_1,b_2$ and $u,v$ are interchanged and $y\rightarrow y^{-1}$. The same is true for the integral $I_2$ if it is evaluated over any of the sets $\{\Delta >0\}$ or $\{\Delta <0\}$.\\[0.5em]
 If $\Delta >0$, we can rewrite $E_2$ as 
 $$E_2=A\cdot (y^\alpha +y^{-\alpha} )+B\cdot y^\alpha ,$$
 where $A>0$ and $B\geq 0$ only depend on $u,v,b_1$ and $b_2$. It follows that $E_2$, and thereby also $E_1$, is a Bernstein function wrt $y^\alpha +y^{-\alpha}$ in this case.\\[0.5em] In the opposite case $\Delta <0$ we get $E_2=
 A\cdot (y^\alpha +y^{-\alpha} )+B\cdot y^{-\alpha }$. It follows from the substitutions above that $I_2$ is unchanged wrt $y\rightarrow y^{-1}$ and we are back in the first case. We conclude that $E_1$ is a Bernstein function wrt $y^\alpha +y^{-\alpha}$ and $I_2$ is CM wrt $ y^\alpha +y^{-\alpha}$.\\[0.5em]
 Then $I_2$ is also CM wrt $y+y^{-1}$ by Bondesson \cite{B1}, Ex. 4.3.4, p. 69, since $y^\alpha+y^{-\alpha}$ is a Bernstein function wrt $y+y^{-1}$. By Bernstein's Theorem, $I_2$ can be represented by a Laplace transform $I_2= \int\limits _0 ^\infty  e^{-\lambda \cdot (y+y^{-1})}d\nu(\lambda )$, for a non-negative Borel measure ${\nu }$. Inserting this formula into $H_2$ then gives
 $$ H_2\sim \int\limits _0^\infty\int\limits _0^\infty e^{- \big( sx^q\cdot (t\cdot y+t^{-1}\cdot y^{-1})+\lambda\cdot (y+y^{-1} )\big) }\, \frac{dy}{y}\, d\nu (\lambda ).$$
 The exponent (with reversed sign) is a linear combination of $y$ and $y^{-1}$ and equals
 $$y\cdot (sx^q\cdot t+\lambda )+y^{-1}\cdot (sx^q\cdot t^{-1}+\lambda )=\rho\cdot \big (s^2x^{2q}+\lambda^2 +sx^q\lambda \cdot (t+t^{-1})\big),$$
 after the substitution putting the second term equal to $1/\rho$. This proves that $H_2$ is CM wrt $t+t^{-1}$ and 
 then $(X_1+X_2)^q\sim GGC$ by Proposition 1, which completes the proof of Theorem 1 in the case $n=2$.\\[0.5em]
Let $n\geq 3$ be an arbitrary integer and let $S_n$ be the sum of the $n$ independent gamma variables defined in (2). We start from the Laplace transform $\phi_{S_n^q}$ of $S_n^q$, 
$$\phi_{S_n^q }(x)\sim \int\limits_0^\infty  e^{-sx^q}\cdot x^{n\beta -1}\cdot    \int\limits _0^1\cdots \int\limits _0^1 e^{-x\cdot C_n({\bf u})}\cdot B_n({\bf u})\, du_1\cdots du_{n-1}\, dx.$$
In analogy with the case $n=2$, we define $H_n=\phi_{S_n^q}(st)\cdot \phi_{S_n^q}(\frac{s}{t})$ as a product of two such integrals and get in analogy with (5)
$$H_n\sim \int\limits_0^\infty\int\limits_0^\infty e^{-\big(sx^q\cdot (t\cdot y+t^{-1}\cdot y^{-1} )\big)} \cdot  x^{2n\beta -1}\cdot$$
 $$\int\limits_0^1\dots \int\limits_0^1 e^{-x\cdot \big (C_n({\bf u})y^\alpha+C_n({\bf v})y^{-\alpha}\big )}\cdot B_n({\bf u})\cdot B_n({\bf v})\, du_1\cdots du_{n-1}dv_1\cdots dv_{n-1}\, dx\frac{dy }{y },$$
 after a hyperbolic change of variables $x\rightarrow x\cdot y$, $y\rightarrow x/y$ and a substitution $y\rightarrow y^\alpha$. We denote the inner integral in $H_n$ by $I_n$. \\[0.5em]
 Now we assume that $I_{n-1}$ is CM wrt $y^\alpha+y^{-\alpha}$ for any sum of $n-1$ independent gamma variables.
 The last two integrals in $I_n$ are equal to
 $$J_n=\int\limits_0^1\int\limits_0^1 e^{-x\cdot E_n} 
 \cdot \big( (1-u_{n-1})\cdot (1-v_{n-1})\big)^{\beta -1} \big( u_{n-1}\cdot  v_{n-1}\big)^{(n-1)\beta -1}\, du_{n-1}\, dv_{n-1},$$ where by (3)
 $$E_n=b_{n-1}\cdot \big( (1-u_{n-1})\cdot Ay^\alpha + (1-v_{n-1})\cdot By^{-\alpha}\big)+$$
$$+b_n\cdot \big( u_{n-1}Ay^\alpha+v_{n-1}By^{-\alpha}\big)$$
 and
$$A=u_1u_2\cdots u_{n-2}, \, \, B=v_1v_2\cdots v_{n-2}.$$
   We apply Lemma 2 with these values on $A$ and $B$. Then $J_n$ is CM wrt $Ay^\alpha +By^{-\alpha}$ and can be represented by a Laplace transform 
 \begin{displaymath}J_n=\int\limits _0^\infty e^{-\lambda\cdot \big(u_1u_2\cdots u_{n-2}\cdot y^\alpha + v_1v_2\cdots v_{n-2}\cdot y^{-\alpha}\big)}\, \nu (d\lambda ),\end{displaymath}
 for some nonnegative Borel measure $\nu$. Now we insert $J_n$  back into $I_n$. Then for every fixed $\lambda >0$, $I_n$ corresponds to $I_{n-1}$ for a sum of $n-1$ independent gamma variables and we conclude that $I_n$ is CM wrt $y^{\alpha}+y^{-\alpha}$, by the induction hypothesis. Recalling that $y^\alpha+y^{-\alpha}$ is a Bernstein function wrt $y +y^{-1}$ and a substitution similar to the one used in the proof for the case $n=2$ then proves that $H_n$ is CM wrt $t+t^{-1}$. We conclude that $S_n^q \sim GGC$ and the proof of Theorem 1 is complete by Proposition 1. 
 \hfill $\Box$\\[1em]
 {\it Proof of Lemma 2.} The proof follows the case $n=2$ in the proof of Theorem 1, with $y^\alpha ,y^{-\alpha}$ replaced by $Ay^\alpha ,By^{-\alpha} $. For the readers convenience we sketch the proof. By (4) we must show that $L =\int\limits_0^1\int\limits_0^1 e^{-  \cdot E}\cdot g(u,v)\, du\, dv, \,$ where
 $$E= b_1\big((1-u)Ay^\alpha +(1-v)By^{-\alpha}\big)+b_2(uAy^\alpha +vBy^{-\alpha}),$$
  is CM wrt $Ay^\alpha +By^{-\alpha}$. We start with substitutions  $u\rightarrow \frac{1}{1+u}$ and $v\rightarrow \frac{1}{1+v}$ to get $L=\int\int e^{-E_1} g(u,v) dudv$, where
 $$E_1=x\cdot  \frac{b_1( u  Ay^\alpha +v  By^{-\alpha})+b_2( v  Ay^\alpha +u  By^{-\alpha})+(b_1uv+b_2) (Ay^\alpha +By^{-\alpha}) }{(1+u)(1+v)}$$ 
  and $g(u,v)=(uv)^{\beta_1-1}\cdot \big((1+u)(1+v)\big)^{-\beta_1-\beta_2}$ and the integration is over $(0,\infty )\times (0,\infty )$.   \\[0.5em]
  It is easy to see that $L$ is unchanged if  $(b_1,b_2)$ and $(u,v)$ are interchanged and $Ay^\alpha \rightarrow By^{-\alpha}$. Without loss of generality we assume that $Ay^\alpha >By^{-\alpha}$ and put $Ay^\alpha +By^{-\alpha}=2s$. Then we get 
  $$Ay^\alpha =s+\sqrt{s^2-AB}\, \textrm{ and }\, By^{-\alpha}=s-\sqrt{s^2-AB}$$
   and note that also $Ay^\alpha$ is a Bernstein function wrt $Ay^\alpha +By^{-\alpha}$. \\[0.5em]
  Let $E_2$ denote the first two terms in the nominator of $E_1$, then we can rewrite $E_2$ as
   $$E_2=(b_1u+b_2v)\cdot Ay^\alpha + (b_1v+b_2u)\cdot By^{-\alpha}$$
   and $\Delta = (b_1u+b_2v)-(b_1v+b_2u)=(b_1-b_2)(u-v)$. The rest of the proof is the same as in the case $n=2$ and is left to the reader. We conclude that $E_1$ is a Bernstein function wrt $Ay^\alpha +By^{-\alpha}$ and $L$ is CM wrt $Ay^\alpha +By^{-\alpha}$, which completes the proof of Lemma 2. \hfill $\Box$
\\[1em]{\it Remark 1.} The old conjecture that $X\sim GGC$ implies $X^q\sim GGC$, $q>1$, is a natural structural property of $GGC$ mentioned in Bondesson \cite{B1}, p. 97. A
different approach was made in Bondesson \cite{B2}, where the Laplace transform of $S_n^q$ is expressed using the product of the densities of the individual rvs  $X_i$, $1\leq i\leq n$. The following sufficient condition for Theorem 1 to be true is given in Bondesson \cite{B2}, Conjecture 2.\\[0.5em] {\it For every $q\geq 1$, $\alpha= 1/q$ and fixed positive numbers $u_1,u_2,\cdots , u_n$ and $\lambda_1,\lambda_2,\cdots , \lambda_n$, the function
$$\int \frac{1}{v_1v_2\cdots v_n}e^{-E}d{\bf v}, \textrm{where}\, E=(\sum u_iv_i^\alpha )^q+(\sum u_iv_i^{-\alpha})^q+\sum \lambda_i(\frac{t}{v_i}+\frac{v_i}{t})$$
is CM wrt $t+t^{-1}$.} 
 \\[0.5em]
{\it Remark 2.} The advantages with the
 method used here compared to the one in Remark 1 is that $I_n$ is inductively defined, the exponent in the integrand of $I_n$ is a linear function wrt $y^\alpha$ and $y^{-\alpha}$, the inner integral $J_n$ has only two variables for all $n$ and  that the method of successive substitutions works here.\\[0.5em]
 {\it Remark 3.} For $n=2$ the integral $I_2$, with $c=1$ can be expressed as a product of two Modified Bessel functions of first order using computer algebra to  be $$I_2=\pi\cdot \Gamma (\beta )^2 \cdot e^{-(y+\frac{1}{y})}\cdot \textrm{BesselI}(\beta ,y/2)\cdot \textrm{BesselI}(\beta ,1/2y). $$ 
 A bold but natural suggestion is that Theorem 1 holds for more general compositions $f\circ X$, $X\sim GGC$, where $f$ belongs to some class of smooth, increasing and convex functions defined on $ [0,\infty)$ and satisfying $f(0)=0$. \section{Applications}
The class $GGC$ is closed wrt sums and products of independent rvs and now also wrt $q-$th powers, $q>1$. This gives the 
following result. 
\begin{theorem} Let $\{ X_i\}_1^n$ be independent rvs, $X_i\sim GGC$, and let $q_i\geq 1$, $1\leq i\leq n$. Then
$\prod\limits _1^n X_i^{q_i}\sim GGC \textrm{ and }\sum\limits _1^n X_i^{q_i}\sim GGC.
$
\end{theorem}
The class of extended generalized gamma convolutions  ($EGGC$) was introduced by Thorin \cite{T3} and consists of limit distributions for sums of independent positive and negative gamma rvs. The symmetric distributions in $ EGGC$ are denoted by $symEGGC$ and are characterized by $X\sim symEGGC$ if and only if  $X=\sqrt{Y}\cdot Z$, for some $Y\sim GGC$ and independent $Z\sim N(0,1)$, see Bondesson \cite{B2}, Ch. 5 or Steutel, van Harn \cite{SH}, Ch. VI, \S 11. Theorem 1 implies that, if $0<\alpha <1$, every $Y\sim GGC$ can be written $Y=Z^\alpha$, for some $Z\sim GGC$. This gives the following extension of Bondesson \cite{B2}, Theorem 2.
\begin{theorem}Let $0<\alpha\leq 2\leq \beta$, then\\ (a) If $Y\sim GGC$ and  $Z\sim N(0,1)$ are independent, then $Y^{1/\alpha} \cdot Z\sim symEGGC$.\\
(b) If $X\sim symEGGC$ there exist $Y\sim GGC$ and an independent $Z\sim N(0,1)$ such that $X=Y^{1/\beta} \cdot Z$.
\end{theorem}
 If $\mathcal A$ and $\mathcal B$ are two classes of distributions, we let $\mathcal A\times \mathcal B$ denote the class of products $X\cdot Y$ of independent rvs $X\sim \mathcal A$ and $Y\sim \mathcal B$. Then we can express Theorem 3 as
 $$GGC^{1/\alpha} \times N(0,1)\subseteq symEGGC\subseteq GGC^{1/\beta}\times N(0,1),\,  $$  $0<\alpha\leq 2\leq \beta $, with equality for $\alpha =\beta =2$.\\[0.5em]
We finally give a new proof of Bondesson \cite{B2} Theorem 3, see the comment on p. 1075.  
\begin{theorem} If $X\sim GGC$ has left extremity $a\geq 0$, then $e^X-e^a\sim GGC$.\end{theorem}
{\it Proof.} If $a=0$, $X\sim GGC$ and $0<r<1$, then $(1+rX)^{\frac{1}{r}}\sim GGC$ by Theorem 1 and
$$Pr((1+rX)^{\frac{1}{r}}\leq 1+u)=Pr(X\leq  [(1+u)^r-1]/r)\rightarrow Pr(e^X\leq 1+u),$$
as $r\rightarrow 0,$ by L'Hopital's rule. Hence $e^X-1\sim GGC$, since $GGC$ is closed wrt weak limits. If $a>0$ we have $e^X-e^a=e^a(e^{X-a}-1)$, $(X-a)\sim GGC$ has left extremity zero and $e^X-e^a\sim GGC$ follows from the first case.\hfill $\Box$ 
 
\section{Acknowledgement} The author thanks Professor Emeritus L. Bondesson for introducing me to the field of Generalized Gamma Convolutions (GGC), the problems treated in this paper and many encouraging discussions and comments on my work. 
  
 \end{document}